\documentclass[11pt, a4paper]{article}
\usepackage{amsmath, amsthm, amssymb, eucal}

\makeatletter
\renewcommand\section{\@startsection{section}{1}{\z@}
 						{-3.5ex \@plus -1ex \@minus -.2ex}
						{2ex \@plus.2ex}
						{\large\bfseries}}
\renewcommand\subsection{\setcounter{subsection}{\value{equation}}
						\stepcounter{equation}
					\@startsection{subsection}{2}{\z@}
                          {1.75ex \@plus.5ex \@minus.2ex}
                           {-.4em}
                           {\textit}}
\def\@seccntformat#1{\@ifundefined{#1@cntformat}
	{\csname the#1\endcsname\quad} 
	{\csname #1@cntformat\endcsname}} 
\def\section@cntformat{\thesection.~} 
\def\subsection@cntformat{(\thesubsection)\ }
\makeatother

\numberwithin{equation}{section}
\theoremstyle{plain}
\swapnumbers
\newtheorem{theorem}[equation]{Theorem}
\newtheorem{conjecture}[equation]{Conjecture}
\newtheorem{amplification}[equation]{Amplification}
\newtheorem{proposition}[equation]{Proposition}
\theoremstyle{definition}
\newtheorem{definition}[equation]{Definition}
\theoremstyle{remark}
\newtheorem{remark}[equation]{Remark}

\newcommand{\cA}{\mathcal{A}}
\newcommand{\cE}{\mathcal{E}}
\newcommand{\cG}{\mathcal{G}}
\newcommand{\cO}{\mathcal{O}}
\newcommand{\cP}{\mathcal{P}}
\newcommand{\frg}{\mathfrak{g}}
\newcommand{\frt}{\mathfrak{t}}
\newcommand{\frM}{\mathfrak{M}}
\newcommand{\bC}{\mathbb{C}}
\newcommand{\bQ}{\mathbb{Q}}
\newcommand{\bZ}{\mathbb{Z}}
\newcommand{\xt}{{}^\tau\!}
										
\begin{document}

\title{\textbf{$K$-theory of the moduli of principal bundles on a surface 
and deformations of the Verlinde algebra}}
\author{Constantin Teleman\\ \small DPMMS, Wilberforce Road,
	\\ \small Cambridge, CB3 0WB, UK\\ \small teleman@dpmms.cam.ac.uk
} 
\maketitle

\begin{quote}
\abstract{\noindent I conjecture\footnote{This note is based on my talk 
at the Segalfest. Recent discussions with Chris Woodward have been 
pointing toward a proof of the main statements. However, I prefer to 
err on the side of caution and preserve their conjectural status.} 
that index formulas for $K$-theory classes on the moduli of 
holomorphic $G$-bundles over a compact Riemann surface $\Sigma$ are 
controlled, in a precise way, by Frobenius algebra deformations 
of the Verlinde algebra of $G$. The Frobenius algebras in question are 
twisted $K$-theories of $G$, equivariant under the conjugation action, 
and the controlling device is the equivariant Gysin map along the 
``product of commutators" from $G^{2g}$ to $G$. The conjecture is 
compatible with naive virtual localization of holomorphic bundles, 
from $G$ to its maximal torus; this follows by localization in twisted 
$K$-theory.}
\end{quote}

\section{Introduction} 
Let $G$ be a compact Lie group and let $M$ be the moduli space of flat 
$G$-bundles on a closed Riemann surface $\Sigma$ of genus $g$. By 
well-known results of Narasimhan, Seshadri and Ramanathan \cite{ns}, 
\cite{ram}, this is also the moduli space of stable holomorphic principal 
bundles over $\Sigma$ for the complexified group $G_\bC$; as a complex 
variety, it carries a fundamental class in complex $K$-homology. This 
paper is concerned with index formulas for vector bundles over $M$. 
The analogous problem in cohomology---integration formulas over $M$ 
for top degree polynomials in the tautological generators---has been 
extensively studied \cite{new}, \cite{kir}, \cite{don}, \cite{thad}, 
\cite{wit}, and, for the smooth versions of $M$, the moduli of vector 
bundles of fixed degree co-prime to the rank, it was completely solved 
in \cite{jefkir}. In that situation, the tautological classes generate 
the rational cohomology ring $H^*(M;\bQ)$. Knowledge of the integration 
formula leads to the intersection pairing, and from here, Poincar\'e 
duality determines this ring  as the quotient of the polynomial ring 
in the tautological generators by the null ideal of the pairing. 

For smooth $M$, index formulas result directly from the Riemann-Roch 
theorem and the integration formula. This breaks down in the singular 
case, and no index formula can be so obtained, for groups other than 
$SU(n)$. Follow-up work \cite{kie} has extended the results of \cite
{jefkir} to the study of the duality pairing in intersection cohomology, 
for some of the singular moduli spaces. Whether a useful connection 
to $K$-theory can be made is not known; so, to that extent, the formulas 
I propose here are new. But I should point out the novel features of 
the new approach, even in the smooth case.

More than merely giving numbers, the conjecture posits a structure to
these indices; to wit, they are controlled by (finite-dimensional)
Frobenius algebras, in the way the Verlinde algebras control the index
of powers of the \textit{determinant line bundle} $D$. The Frobenius 
algebras in question are formal deformations of Verlinde algebras. 
This is best explained by the twisted $K$-theory point of view, \cite
{fht1}, \cite{fht2}, which identifies the Verlinde algebra with a twisted 
equivariant $K$-theory $\xt K_G^*(G)$. The determinantal twistings
appearing in that theorem (cf.~\S\ref{13} below) correspond to powers 
of $D$. Other $K$-theory classes, involving index bundles over $\Sigma$  
(see \ref{6}.iii), relate to higher twistings in $K$-theory, and these 
effect infinitesimal deformations of the Verlinde algebra. If ordinary
(=determinantal) twistings can be represented by gerbes \cite{ozgerb}, 
higher twistings are realized by what could be called virtual gerbes, 
which generalize gerbes in the way one-dimensional virtual bundles 
generalize line bundles. In this picture, $K$-theory classes over $M$ 
of virtual dimension one are automorphisms of virtual gerbes, and arise
by comparing two trivializations of a twisting $\tau$ for $K_G(G)$ 
(\S\ref{39}). Relevant examples are the transgressions over $\Sigma$ 
of delooped twistings for $BG$ (\S\ref{21}).\footnote{For the expert, 
we mean $G$-equivariant $B^2BU_\otimes$-classes of a point.} All
these $K$-classes share with $D$ the property of being ``multiplicative 
in a piece of surface". 

The first formulation (\ref{36}) of the conjecture expresses the 
index of such a $K$-class over the moduli of $G$-bundles as the 
\textit{partition function} for the surface $\Sigma$, in the 2D 
topological field theory defined by $\xt K_G^*(G)$. This is a sum 
of powers of the \textit{structure constants} of the Frobenius algebra, 
for which explicit formulas can be given (Thm.~\ref{31}). There is a 
restriction on the allowed twistings, but they are general enough to 
give a satisfactory set of $K$-classes.

We can reformulate the conjecture in (\ref{41}), (\ref{43}) by 
encoding part of the Frobenius algebra structure into the \textit
{product of commutators map} $\Pi:G^{2g}\to G$. This map has the 
virtue of lifting the transgressed twistings for $K_G(G)$ (\S\ref{21})
to trivializable ones, which allows one to identify $K_G^*(G^{2g})$ 
with its $\tau$-twisted version. In particular, we get a class ${}^\tau
 1$ in $\xt K_G^*(G^{2g})$. The conjecture then asserts that the index 
of the $K$-class associated to $\tau$ over the moduli of $G$-bundles 
equals the Frobenius algebra trace of the Gysin push-forward of 
${}^\tau 1$ along $\Pi$.

Having (conjecturally) reduced this index to a map of compact manifolds,
ordinary localization methods allows us to express the answer in terms
of the maximal torus $T$ and Weyl group of $G$. This reduction to $T$, 
it turns out, can be interpreted as a virtual localization theorem from 
the moduli of holomorphic $G_\bC$-bundles to that of holomorphic 
$T_\bC$-bundles. (The word ``virtual" reflects the use of the virtual 
normal bundle, defined by infinitesimal deformations). For this 
interpretation, however, it turns out that we must employ the moduli 
$\frM = \frM_G$, $\frM_T$ of \textit{all} $G_\bC$-and $T_\bC$-bundles, 
not merely the semi-stable ones. These moduli have the structure of 
\textit{smooth stacks}, with an infinite descending stratification by 
smooth algebraic substacks. Even the simplest case, the Verlinde formula, 
cannot be reduced to a single integral over the variety of topologically 
trivial $T_\bC$-bundles; the correct expression arises only upon summing 
over all topological $T$-types. (Recall that the non-trivial $T$-types 
define unstable $G_\bC$-bundles). 

Now is the right time to qualify the advertised statements. The fact 
that the Verlinde formula, the simplest instance of our conjecture, 
expresses the indices of positive powers $D$ over $M$ is a fortunate 
accident. It is an instance of the ``quantization commutes with reduction" 
conjecture of Guillemin and Sternberg \cite{gs}, which in this case 
\cite{quant} equates the indices of positive powers of $D$ over $M$ and 
over the stack $\frM$ of all holomorphic $G_\bC$-bundles. This does not 
hold for more general $K$-theory classes, for which there will be 
contributions from the unstable Atiyah-Bott strata, and our deformed 
Verlinde algebras really control not the index over $M$, but that over 
$\frM$. This incorporates information about the moduli of flat $G$-bundles, 
and the moduli of flat principal bundles of various subgroups of $G$. 
In other words, the index information which assembles to a nice 
structure refers to the stack $\frM$ and not to the space $M$.

Hence, the third formulation of the conjecture (\ref{51}) expresses the 
index of any admissible $K$-class (Def.~\ref{7}) over the moduli stack 
$\frM$ of all holomorphic $G_\bC$-bundles over $\Sigma$ by virtual 
localization to the stack of holomorphic $T_\bC$-bundles. This involves 
integration over the Jacobians, summation over all degrees, leading to a 
distribution on $T$, and finally, integration over $T$ (to extract the 
invariant part). In \S5, these steps are carried out explicitly for 
the group $SU(2)$.

However, even if our interest lies in $M$ (which, our approach
suggests, it should not), all is not lost, because a generalization of
``quantization commutes with reduction" (first proved in \cite{tz}, 
for compact symplectic manifolds) asserts, in this case, the equality 
of indices over $M$ and $\frM$, after a large $D$-power twist, for 
the class of bundles we are considering.\footnote{An explicit bound 
for the power can be given, linear in the highest weight.} This 
follows easily from the methods of \cite{quant}, but there is 
at present no written account. Because $M$ is projective algebraic, 
the index of $\cE\otimes D^{\otimes n}$ over $M$, for any coherent 
sheaf $\cE$, is a polynomial in $n$; so its knowledge for large $n$ 
determines it for all $n$, including, by extrapolation, $n=0$. Thus, 
the information contained in $\frM$, which combines index information 
for the moduli of bundles of subgroups of $G$, can be disassembled 
into its constituent parts; the leading contribution, as $n\to \infty$, 
comes from $M$ itself. When $M$ is smooth, the ``$\cE$-derivatives" 
of the $n\to \infty$ asymptotics of the index of $\cE\otimes D^{\otimes 
n}$ give integration formulas for Chern polynomials of $\cE$ over $M$; 
and the author suggests that the Jeffrey-Kirwan residue formulas for
the integrals can be recovered in this manner. What one definitely 
recovers in the large level limit are Witten's conjectural formulas 
\cite{wit}. Indeed, there is evidence that the relevant field 
theories are topological limits of Yang-Mills theory coupled to the 
WZW model (in other words, the $G/G$ coset model with a Yang-Mills 
term); this would fit in nicely with the physical argument of \cite{wit}.

\section{The moduli space $M$, the moduli stack $\frM$ and 
admissible $K$-classes} 
In this section, we recall some background material; some of it is 
logically needed for the main conjecture, but mostly, it sets the stage 
for my approach to the question. This is anchored in Thm.~\ref{8}.
\subsection{}\label{1} Recall the set-up of \cite{atbott}: let $\cA$ be the 
affine space of smooth connections, and  $\cG$  the group of smooth 
gauge transformations on a fixed smooth principal $G$-bundle $P$ over
$\Sigma$. The $(0,1)$-component of such a connection defines a
$\bar \partial$-operator, hence a complex structure on the principal
$G_\bC$-bundle $P_\bC$ associated to $P$. We can identify $\cA$ with 
the space of smooth connections of type $(0,1)$ on $P_\bC$; the latter 
carries an action of the complexification $\cG_\bC$ of $\cG$, and the 
quotient $\cA\!\left/\cG_\bC\right.$ is the set of isomorphism classes 
of holomorphic principal $G_\bC$ -bundles on $\Sigma$ with underlying 
topological bundle $P_\bC$.
\subsection{}\label{2} The space $\cA$ carries a $\cG_\bC$-equivariant 
stratification, according to the instability type of the holomorphic 
bundle. The semi-stable bundles define the open subset $\cA^0$, whose 
universal Hausdorff quotient by $\cG_\bC$ is a projective algebraic 
variety $M$, the moduli space of semi-stable holomorphic $G_\bC$-bundles 
over $\Sigma$. The complex structure is descended from that of $\cA^0$, 
in the sense that a function on $M$ is holomorphic in an open subset 
if and only if its lift to $\cA$ is so. We can restate this, to avoid 
troubles relating to holomorphy in infinite dimensions. The gauge 
transformations that are based at one point $*\in\Sigma$ act freely 
on $\cA^0$, with quotient a smooth, finite-dimensional algebraic variety 
$M_*$. This is the moduli of semi-stable bundles with a trivialization 
of the fibre over $*$. Its algebro-geometric quotient under the residual 
gauge group $G_\bC$ is $M$. The other strata  $\cA^\xi$ are smooth, locally 
closed complex submanifolds, of finite codimension; they are labeled by 
the non-zero dominant co-weights $\xi$ of $G$, which give the 
destabilizing type of the underlying holomorphic $G_\bC$ -bundle. 
The universal Hausdorff quotient of $\cA^\xi$ can be identified with a 
moduli space of semi-stable principal bundles under the centralizer
$G_{(\xi)}$ of $\xi$ in $G_\bC$.
\subsection{}\label{3} The stack $\frM$ of all holomorphic $G_\bC$-bundles 
over $\Sigma$  is the homotopy quotient $\cA\!\left/\cG_\bC\right.$. As 
such, it seems that we are using new words for an old object, and so it 
would be, were our interest confined to ordinary cohomology, $H^*(\frM;\bZ)$. 
However, we will need to discuss its $K$-theory, and the index map to $\bZ$. 
The abstract setting for this type of question is a homotopy category of 
analytic spaces or algebraic varieties (see e.g. Simpsonês work \cite{sim}, 
and, with reference to $\frM$, \cite{bwb}). Fortunately, little of that 
general abstraction is necessary here. It turns out that  $\frM$ is 
homotopy equivalent to the quotient, by $G$-conjugation, of a principal 
$\Omega G$-bundle over $G^{2g}$. (This is the homotopy fibre of $\Pi: 
G^{2g}\to G$; cf.~\S\ref{38}). This follows from Segal's \textit{double
coset presentation} of $\frM$ (\cite{ps}, Ch.~8). As a result, there is 
a sensible topological definition of $K^*(M)$, which makes it into 
an inverse limit of finite modules for the representation ring
$R_G$ of $G$.
\subsection{}\label{4} As in the cohomological setting of Atiyah-Bott, 
this $K^*(M)$ can be shown to surject onto the $\cG_\bC$-equivariant 
$K$-theory of $\cA^0$. The latter can be defined as the $G$-equivariant 
version of $K^*(M_*)$, for the variety $M_*$ of \S\ref{2}. This is as 
close as we can get to $K^*(M)$. When the action of $G$ on $M_*$ is free, 
the two groups coincide, but this only happens when $G=PU(n)$ and the 
degree of our bundle $P$ is prime to $n$. However, some relation between
$K_G^0(M_*)$ and $M$ always exists. Namely, every holomorphic, $G_\bC
$-equivariant bundle  $\cE$ over $M_*$ has an invariant direct image
$q_*^G(\cE)$, which is a coherent analytic sheaf over $M$. (This is the 
sheaf of $G$-invariant holomorphic sections along the fibres of the 
projection $q:M_*\to M$). The coherent sheaf cohomology groups of 
$q_*^G(\cE)$ are finite dimensional vector spaces, and the alternating 
sum of their dimensions is our definition of the $G$-invariant index 
of $\cE$ over $M_*$.
\subsection{}\label{5} A similar construction, applied to a stratum
$\cA^\xi$, allows us to define the $\cG_\bC$-invariant index of 
holomorphic vector bundles over it. Because $\cA$ is stratified by the
$\cA^\xi$, there is an obvious candidate for the $\cG_\bC$-invariant 
index of holomorphic vector bundles over $\cA$, as a sum over all $\xi$. 
(Contributions from the normal bundle must be taken into account; see
\cite{quant}, \S9). This sum may well be infinite. The task, therefore, 
is to identify a set of admissible $K$-theory classes, for which the 
sum is finite; we then define that sum to be the index over $\frM$. 
For a good 
class of bundles, this can be done, and shown to agree with a more 
abstract global definition, as the \textit{coherent-sheaf cohomology 
Euler characteristic over the algebraic site of the stack $\frM$}. 
I shall not prove any of the assertions above, instead will take a 
low-brow approach and define directly the set of admissible $K$-theory 
classes. When $G$ is simply connected, they turn out to generate 
a dense subring of $K^*(M;\bQ)$, rather in the way that a polynomial 
ring is dense in its power series completion; and their restrictions 
to the semi-stable part generate $K_G^*(M_*)\otimes\bQ$. (This can
be deduced from the cohomological result of \cite{atbott}, by using 
equivariant Chern characters).
\subsection{}\label{6} Note, first, that the pull-back of the bundle
$P_\bC$ to $\Sigma\times\cA$ carries a natural, $\cG_\bC$-equivariant, 
holomorphic structure, as defined by the $(0,1)$-part of the universal 
connection along $\Sigma$. We might call this the universal bundle on 
$\Sigma\times \frM$. A representation $V$ of $G$ defines an associated 
holomorphic, $\cG_\bC$-equivariant vector bundle $E^*V$ on  $\Sigma
\times\cA$. (We think of $E$ as the classifying map of the universal
bundle to $BG$). Let $\pi:\Sigma \times \cA\to \cA$ be the projection, 
and fix a square root $K^{1/2}$ of the canonical bundle on $\Sigma$. 
We now associate the following objects to $V$, which we shall call the
\textit{tautological classes} in $K^*(M)$.
\begin{enumerate}
\item For a point $x\in\Sigma$, the restriction $E_x^*V$ of $E^*V$ to 
$\{ x\}\times\frM$; 
\item The \textit{index bundle} $\alpha(V):= R^*\pi_*(E^*V\otimes 
K^{1/2})$ along  $\Sigma$  over $\cA$; 
\item For any class $C\in K_1(\Sigma)$, its slant product with $E^*V$
(the index of $E^*V$ along a $1$-cycle).
\end{enumerate}
Object (i) is an equivariant holomorphic vector bundle, objects (ii) and
(iii) are equivariant $K^0$ and $K^1$ classes over $\cA$, the misnomer 
``bundle" in (ii) notwithstanding. For example, we can represent 
(ii) by a $\cG_\bC$-equivariant Fredholm complex based on the relative
$\bar\partial$-operator. (The square root of $K$ leads to the Dirac, 
rather than $\bar \partial$ index, and the notation $\alpha(V)$ stems 
from the Atiyah index map, which we get when $\Sigma$ is the sphere.) 
We shall not consider type (iii) objects in this paper, so we refrain 
from analyzing them further. Note that the topological type of the
bundles in (i) is independent of $x$; we shall indeed see that their 
index is so as well.

Any reasonable definition of $K^*(M)$ should include the tautological 
classes, but another distinguished object plays a crucial role:

(iv) The determinant line bundle $D(V):=\det R^*\pi _*(E^*V)$ over $\cA$.

\noindent This is a holomorphic, $\cG_\bC$-equivariant line bundle, 
or a holomorphic line bundle over $\frM$. When $G$ is semi-simple, such 
line bundles turn out to be classified by their Chern classes in $H^2(M;
\bZ)$ \cite{bwb}; in the case of $D(V)$, this is the transgression along 
$\Sigma$ of $c_2(V)\in H^4(BG)$. Not all line bundles are determinants,
but they are fractional powers thereof. The convex hull of the $D(V)$ 
define the semi-positive cone in the group of line bundles; its interior 
is the positive cone. In particular, when $G$ is simple and simply connected,
$H^2(\frM;\bZ)\cong \bZ$, and the positive cone consists of positive powers 
of a single $D$; for $G=SU(n)$, this is $D(\bC^n)$, for the standard 
representation.
\begin{definition}\label{7} An \textit{admissible class} in $K^*(\frM)$ 
is a polynomial in the tautological classes and the semi-positive line 
bundles. 
\end{definition}

For simply connected $G$, an admissible $K$-class is a finite sum of 
terms $p_n\otimes D^{\otimes n}$, where $n\ge 0$ and $p_n$ is a
polynomial in the objects (i)-(iii) above. We can actually allow 
some small\footnote{Larger than the negative of the dual Coxeter number}
negative values of $n$, but the index of such classes turns out to 
vanish, so little is gained. The following theorem allows our 
approach to $K^*(M)$ to get off the ground. To simplify the statements, 
we assume that $G$ is semi-simple.
\begin{theorem}\label{8}
\begin{trivlist} 
\item\textup{(i)} The coherent sheaf cohomology groups, over the
algebraic site of the stack $\frM$, of any admissible class
$\cE\in K^*(\frM)$, are finite-dimensional, and vanish in high degrees. 
\item\textup{(ii)} The index $\mathrm{Ind}(\frM;\cE)$ over $\frM$, defined 
as the alternating sum of cohomology dimensions, is also expressible 
as a sum of index contributions over the Atiyah-Bott strata $\cA^\xi$. 
Each contribution is the index of a coherent sheaf over a moduli 
space of semi-stable $G_{(\xi)}$-bundles. These contributions 
vanish for large $\xi$. 
\item\textup{(iii)} For sufficiently positive $D$, all $\xi\ne 0
$-contributions of $\cE\otimes D$ to the index vanish. 
\item\textup{(iv)} Hence, for sufficiently positive $D$, $\mathrm{Ind}
\left(\frM;\cE\otimes D \right) = \mathrm{Ind}\left(
M;q_*^G(\cE\otimes D) \right)$.
\end{trivlist}
\end{theorem}
\begin{proof} (Sketch) For a product of ``evaluation bundles" (\ref{6}.i), 
the results were proved in \cite{bwb} and \cite{quant}. This generalizes 
immediately to a family of bundles, parametrized by a product of copies of
$\Sigma$, and integration along the curves leads to index bundles 
(\ref{6}.ii). Slant products with odd $K$-homology classes on that 
product of Riemann surfaces lead to the same conclusion for arbitrary 
admissible $\cE$.
\end{proof}
\begin{remark}\label{9}
As noted in the introduction, the theorem allows us to determine the 
index of $q_*^G(\cE)$ over $M$, for any admissible $\cE$, if index
formulas over $\frM$ are known. Indeed, suitable positive line bundles 
$D$ descend to $M$, and so $\mathrm{Ind}\left(M ; q_*^G(\cE\otimes D^
{\otimes n})\right)$ is a polynomial in $n$. Its value at $n=0$ can then 
be determined from large $n$, where it agrees with $\mathrm{Ind}\left(
\frM;\cE\otimes D^{\otimes n}\right)$.
\end{remark}

\section{Twistings and higher twistings in $K$-theory} 
We start with some background on twisted $K$-theory and its equivariant 
versions. The statements are of the ``known to experts" kind, but 
unfortunately, references do not always exist. They will be proved 
elsewhere.
\subsection{}\label{10} Let $X$ be a compact, connected space. Units in 
the ring $K^*(X)$, under tensor product, are represented by the virtual 
vector bundles of dimension $\pm1$ . A distinguished set of units in the 
1-dimensional part $GL_1^+$ is the Picard group $\mathrm{Pic}(X)$ of topological 
line bundles; it is isomorphic to $H^2(X;\bZ)$, by the Chern class. In 
the other direction, the determinant defines a splitting
\begin{equation}\label{11}			
GL_1^+\left(K^*(X) \right) \cong \mathrm{Pic}(X)\times 
	SL_1\left(K^*(X)\right), 
\end{equation}
\noindent where the last factor denotes 1-dimensional virtual
bundles with trivialized determinant line.

We shall ignore here the twistings  coming from the group $\{\pm 1\}$.
The splitting (\ref{11}) refines to a decomposition of the spectrum
$BU_\otimes$ of 1-dimensional units in the classifying spectrum for 
complex $K$-theory \cite{mst}; in self-explanatory notation, we have a 
factorization
\begin{equation}\label{12}			
BU_\otimes \cong K\left(\bZ;2\right)\times BSU_\otimes.
\end{equation}
\subsection{}\label{13} A twisting of complex $K$-theory over $X$ 
is a principal $BU_\otimes$-bundle over that space. By (\ref{12}), 
this is a pair $\tau =\left( {\delta ,\chi } \right)$ consisting of 
a \textit{determinantal twisting} $\delta$, which is a $K(\bZ;2)$-principal 
bundle over $X$, and a \textit{higher twisting} $\chi$, which is a
$BSU_\otimes$-torsor. Twistings are classified, up to isomorphism, by
a pair of classes $[\delta ]\in H^3(X;\bZ)$ and $[\chi]$ in the
generalized cohomology group $H^1(X;BSU_\otimes)$. This last
group has some subtle features over $\bZ$; rationally, however, 
$BSU_\otimes$ is a topological abelian group, isomorphic to
$\prod_{n\ge 2}K(\bQ;2n)$ via the logarithm of the Chern character 
$ch$. We obtain the following.
\begin{proposition}\label{14}
Twistings of rational $K$-theory over $X$ are classified, up to 
isomorphism, by the group $\prod_{n>1}H^{2n+1}(X;\bQ)$. \qed
\end{proposition}
\begin{remark}\label{15} 
The usual caveat applies: if $X$ is not a finite complex, we
rationalize the coefficients before computing cohomologies.
\end{remark}

The twistings in Prop.~\ref{14}, of course, are also the twistings 
for rational cohomology with coefficients in formal Laurent series
$\bQ((\beta))$ in the Bott element $\beta$ of degree $(-2)$. This is
not surprising, as the classifying spectra $BU\otimes\bQ$ and $K\left(
\bQ((\beta )),0\right)$ for the two theories are equivalent 
under $ch$. The isomorphism extends naturally to the twisted theories, 
by a twisted version of the Chern character, as in \S3 of \cite{fht1}, 
where determinantal twistings were considered:
\begin{proposition}\label{16}
There is a natural isomorphism $\xt ch:\xt K^*\left(X;\bQ\right)
\to \xt H^*\left(X;\bQ((\beta))\right)$. \qed
\end{proposition}
\begin{remark}\label{17}
The strength of the proposition stems from computability of the 
right-hand side. Let $(A^\bullet ,d)$ be a DGA model for the rational 
homotopy of $X$, $\eta =\eta ^{(3)}+\eta ^{(5)}+\ldots$ a cocycle 
representing the twisting, decomposed into graded parts. If we define
$\eta':=\beta \eta ^{(3)}+\beta ^2\eta ^{(5)}+\ldots$, then it turns 
out that $\xt H^*\left(X;\bQ((\beta ))\right)$ is the cohomology of 
$A^\bullet((\beta ))$ with modified differential $d+\eta'\wedge$. 
The latter can be computed by a spectral sequence, commencing at 
$E_2$ with the ordinary $H^*\left(X;\bQ((\beta))\right)$, and with 
third differential $\beta \eta^{(3)}$ (cf.~\cite
{fht1}).
\end{remark}
\subsection{}\label{18} 
Thus far, the splitting (\ref{12}) has not played a conspicuous role:
rationally, all twistings can be treated uniformly, with $\log ch$ 
playing the role that $c_1$ plays for determinantal twistings. Things 
stand differently in the equivariant world, when a compact group $G$
acts on $X$. There are equivariant counterparts to (\ref{10}) and (\ref
{12}), namely, a spectrum $BU_\otimes ^G$ of equivariant $K$-theory units, 
factoring into the equivariant versions of $K(\bZ;2)$ and $BSU_\otimes$. 
However, the group analogous to (\ref{14}), $\prod _{n>1}H_G^{2n+1}
\left(X;\bQ\right)$, no longer classifies twistings for rational
equivariant $K$-theory, but only for its augmentation completion. 
The reason is the comparative dearth of units in the representation 
ring $R_G$ of $G$, versus its completion.

The easily salvaged part in (\ref{11}) is the equivariant Picard group.
Realizing twistings in $H_G^3(X;\bZ)$ by equivariant projective Hilbert
bundles allows a construction of the associated twisted $K$-theory by
Fredholm operators; see \cite{icm} for a relevant example.\footnote{Even 
here, we meet a new phenomenon, in that \textit{integral} twistings are 
required to define the rational equivariant $K$-theory; for instance, the 
torsion part of a twisting in $H_G^3(X;\bZ)$ affects the rational answer.}
One method to include more units (hence more twistings) is to localize 
in the representation ring. This is not an option here; one of our 
operations for the index formula, the trace map, will be integration 
over $G$. Instead, we introduce the extra units by adjoining formal 
variables to $R_G$. Thus, we consider the formal power series ring
$R_G[[t]]$, and the associated $K$-theory $K_G^*(X)[[t]]$ of formal power 
series in $t$, with equivariant vector bundle coefficients. Any series
$\sum t^nV_n$, where $V_0$ is a 1-dimensional representation of $G$, 
is now a unit. This realizes the $G$-spectrum
\begin{equation}\label{19}			
K^G(\bZ;2)\times \left(1+tBU^G[[t]]\right)_\otimes 
\end{equation}
\noindent within the $G$-spectrum of units in $BU^G[[t]]$.
\begin{definition}\label{20}
An \textit{admissible twisting} for $K_G(X)[[t]]$-theory is a torsor 
under the spectrum (\ref{19}) over the $G$-space $X$.
\end{definition}

\noindent Admissible twistings $\tau(t) = (\delta,\chi(t))$ are classified, 
up to isomorphism, by a pair of classes $[\delta ]\in H_G^3(X;\bZ)$ 
and a ``higher" class $[\chi(t)]$ in the generalized equivariant 
cohomology group $\left[X; B\left(1+tBU[[t]]\right)_\otimes\right]^G$.

\subsection{} \label{21} One way to define a higher admissible twisting
over $G$, equivariant for the conjugation action, uses a (twice deloopable) 
\textit{exponential morphism} from $BU^G[[t]]$ to $1 + tBU^G[[t]]$ 
(taking sums to tensor products). Segal's theory of $\Gamma$-spaces 
\cite{seg} shows that sufficiently natural exponential operations on 
the coefficient ring $R_G[[t]]$ define such morphisms. An example is 
the \textit{total symmetric power}
\begin{equation}\label{22}			
V(t)\mapsto S_t\left[ {t\cdot V(t)} \right]=\sum_{n\ge 0}\,t^n\cdot S^n
\left[{V(t)} \right].
\end{equation}
\noindent Allowing rational coefficients, the naive exponential,
\begin{equation}\label{23}			
V(t)\mapsto \exp\left[t\cdot V(t)\right] = \sum_{n\ge 0} 
{t^n/n!\cdot V(t)^{\otimes n}}, 
\end{equation}
\noindent is more closely related to the previous discussion, in the 
sense that completing at the augmentation ideal takes us to the ordinary 
$K$-theory of $BG$, and applying $ch$ to the right leads to the earlier 
identification (\ref{13})--(\ref{14}) of $BSU_\otimes \otimes\bQ$ with 
$\prod _{n\ge2} K(\bQ;2n)$. Whichever exponential morphism we choose, 
Bott periodicity permits us to regard $BU^G[[t]]_\oplus$-classes of a 
point as $B^2BU^G[[t]]_\oplus$-classes, and delooping our morphism 
produces classes in $B^2(1 + tBU^G[[t]])$ of a point. Transgressing 
once gives higher twistings for $K_G(G)$. In this paper, we shall pursue 
the exponential (\ref{23}) in more detail.

\section{The Index from Verlinde algebras}
\subsection{The Verlinde algebra and its deformations.} 
Call a twisting $\tau(t) = (\tau_0,\chi (t))$ \textit{non-degenerate} 
if the invariant bilinear form $h$ it defines on $\frg$, via the 
restriction of $[\tau_0]\in H_G^3(G)$ to $H_T^2\otimes H^1(T)$, is so; 
call it \textit{positive} if this same form is symmetric and positive 
definite. Integrality of $\tau_0$ implies that $h$ defines an isogeny from 
$T$ to the Langlands dual torus, with kernel a finite, Weyl-invariant 
subgroup $F\subset T$. Recall the following result from \cite{fht1}, 
\cite{fht2}, referring to determinantal twisting $\tau = (\tau_0,0)$. 
For simplicity, we restrict to the simply connected case. Let $\sigma$ 
be the twisting coming from the projective cocycle of the Spin 
representation of the loop group; it restricts to the dual Coxeter 
number on $H^3$ for each simple factor.
\begin{theorem}\label{24}
Let $G$ be simply connected. For a positive determinantal twisting $\tau$, 
the twisted K-theory $\xt K_G^{\dim G}(G)$ is isomorphic to the Verlinde
algebra $V_G(\tau -\sigma)$ of $G$, at a shifted level $\tau-\sigma$. 
It has the structure of an integral Frobenius algebra; as a ring, it 
is the quotient of $R_G$ by the ideal of representations whose characters
vanish at the regular points of $F$. The trace form ${}^\tau Tr:R_G\to\bZ$ 
sends $V\in R_G$ to
\begin{equation}\label{25}			
{}^\tau Tr(V) = \sum_{f\in F^{reg}/W} {ch_V(f)\cdot
\frac{\left|\Delta(f)\right|^2}{|F|}}, 
\end{equation}
\noindent where $ch_V$ is the character of $V$, $\Delta (f)$ is the 
Weyl denominator and $\left| F\right|$ is the order of $F$. \qed
\end{theorem}

\begin{remark}\label{26}
\begin{trivlist}\itemsep0ex
\item(i) The $R_G\otimes\bC$-algebra $\xt K_G^{\dim G}(G)\otimes\bC$ 
is supported at the regular conjugacy classes of $G$ which meet $F$, 
and has one-dimensional fibres. 
\item(ii) The trace form (\ref{25}) determines the Frobenius 
algebra: by non-degeneracy, the kernel of the homomorphism from
$R_G$ is the null subspace under the bilinear form $(V,W)\mapsto {}^\tau
Tr(V\otimes W)$. 
\item(iii) After complexifying $R_G$, we can represent Tr by 
integration against an invariant distribution on $G$. The latter is 
the sum of $\delta$-functions on the conjugacy classes in (ii), divided 
by the order of $F$; the factor $\left| {\Delta (f)} \right|^2$ is the 
volume of the conjugacy class. 
\item(iv) The result holds for connected groups with torsion-free $\pi_1$, 
although some care must be taken with the ring structure when the adjoint 
representation does not spin \cite{fht1}. When $\pi_1$ has torsion or 
$G$ is disconnected, $\xt K_G^{\dim G}(G)$ is still the Verlinde algebra, 
but it is larger than the quotient of $R_G$ described in Thm.~\ref{24}; 
so the trace form on $R_G$ no longer determines $V_G$.
\end{trivlist}
\end{remark}

The complexified form of Thm.~\ref{24} was given a direct proof in
\cite{fht1}, using the Chern character to compute the twisted $K$-theory. 
The trace form was introduced ad hoc, using our knowledge of the Verlinde
algebra (\S7 of loc.~cit.). There is in fact no choice on the matter, and
the entire Frobenius structure is determined topologically 
(see the proof of Thm. \ref{31} below). 

We now incorporate higher twistings in Theorem \ref{24}. We tensor with
$\bC$ for convenience.
\begin{theorem}\label{27}
Let G be simply connected, $\tau(t) = (\tau_0,\chi(t))$ an admissible 
twisting with positive determinantal part $\tau_0$. The twisted $K$-theory 
$\xt K_G^{\dim G}(G)\otimes [[t]]$ is a Frobenius algebra, which is a 
quotient of $R_G[[t]]\otimes\bC$, and a flat deformation of the Verlinde 
algebra at level $\tau_0$.
\end{theorem}
\begin{remark}\label{28}
The use of complex higher twistings forces us to tensor with $\bC$ . 
The use of integral twistings of the type (\ref{22}) would lead to a 
similar result over $\bZ$, but as our goal here is an index formula, 
nothing is lost over $\bC$.
\end{remark}
\begin{proof} (Idea) All statements follow by computing the Chern 
character, exactly as in \cite{fht1}; the completions of $\xt K_G^
{\dim G}(G)\otimes \bC[[t]]$ at conjugacy classes in $G_\bC$ are 
calculable by spectral sequences as in (\ref{17}). Away from $F$, 
the $E_2$ term of this sequence is nil. At singular points of $F$, 
this is not so, but the third differential, which stems from the 
determinant of the twisting, is exact; so the limit is null again.
At regular points of $F$, the same third differential resolves one
copy of $\bC[[t]]$ in degrees of the parity $\dim G$, and zero 
otherwise; so the sequence collapses there, and the  abutment 
is a free $\bC[[t]] $-module of rank 1.
\end{proof}
\subsection{}\label{29} As explained in (\ref{26}.ii), the Frobenius 
algebra is completely determined by the trace form ${}^\tau Tr:R_G[[t]]
\to \bC[[t]]$. At $t=0$ , this is given by an invariant distribution 
$\varphi_0$ on $G$, specifically, $1/|F|$ times the sum of $\delta
$-functions on the regular conjugacy classes meeting $F$. We must 
describe how this varies with $t$. Recall that the determinantal part 
$\tau(0) = (\tau_0,0)$ of the twisting defines an invariant metric $h$ 
on $\frg$. We will now associate to the unit $\exp(t\cdot V)$ a 
one-parameter family of conjugation-invariant coordinate changes on 
the group $G$. More precisely, this is a (formal) path in the complexified 
(formal) group of automorphisms of the variety $G/G$; or, even more 
precisely, a formal 1-parameter family of automorphisms of the 
representation ring $R_G\otimes\bC$. Because $G/G = T/W$, it suffices 
to describe this on the maximal torus $T$. In flat coordinates 
$\exp(\xi)$, where $\xi \in \frt$, this is
\begin{equation}\label{30}				
\xi \mapsto \xi +t\cdot
\nabla \left[ch_V\left(\exp(\xi)\right)\right], 
\end{equation}
\noindent the gradient being computed with respect to the metric $h$ 
on $\frt$.
\begin{theorem}\label{31}
The trace form ${}^\tau Tr:R_G[[t]]\mapsto \bC[[t]]$ is integration 
against the invariant distribution $\varphi_t$ on $G_\bC$, obtained 
from $\varphi_0$, the distribution associated to $\tau_0=(\tau_0,0)$, 
by the formal family of coordinate changes (\ref{30}) on $G_\bC$.
\end{theorem}
\begin{remark}\label{32}
It is no more difficult to give the formula for more general exponential
morphisms $V\mapsto \Phi_t(V)$. We assume $\Phi_t$ compatible with 
the splitting principle, via restriction to the maximal torus: in other 
words, its value on any line $L$ is a formal power series in $t$, with 
coefficients Laurent polynomials in $L$. Then, $\log \Phi_t$ extends 
by linearity to an additive map $\frt\otimes R_T\to R_T[[t]]$, and the 
required change of coordinates is $\xi \mapsto \xi +\log \Phi_t
\left[\nabla ch_V\right]$, the metric being used to define the gradient.
For instance, the symmetric power twisting (\ref{22}) arises from 
$\Phi_t(L) = (1-tL)^{-1}$; when $G = S^1$, $ch_V(u) = \sum c_nu^n$,
and we take level $h$ for the determinantal part, we get the change of 
variable $u \mapsto u\cdot \prod (1 - tu^n)^{-nc_n/h}$.
\end{remark}
\begin{proof}(Sketch) The 2D field theory structure of $\xt K_G(G)$
requires the trace form to be the inverse of the bilinear form which
is the image of $1 \in K_G(G)$ in $\xt K_G(G)^{\otimes2}$, under the
anti-diagonal morphism of spaces $G \to G\times G$ (and the diagonal
inclusion of the acting groups). By localization to the maximal torus
$T$, it suffices to check the proposition for tori. (The Euler class
of the inclusion $T \subset G$ is responsible for the factor $|\Delta|^2$). 
Now, the twisting $\tau$ enters the computation of this direct image 
only via the holonomy representation $\pi_1(T) \to GL_1(R_T[[t]])$ it 
defines. Via the metric $h$, the determinantal twisting $\tau_0$ assigns 
to any $p \in\pi_1(T)$ a weight of $T$, which gives a unit in $R_T$, 
and this defines the holonomy representation for $\tau_0$. The change 
of coordinates (\ref{30}) has the precise effect of converting this 
holonomy representation to the one associated to $\tau$. 
\end{proof}
\subsection{Index formulas.}\label{33} In the simply connected case, 
the class $[\tau _0]\in H_G^3(G;\bZ)$ determines a unique holomorphic 
line bundle over $\frM$ (cf.~\ref{6}), which we call $\cO(\tau_0)$. 
The simplest relation between the Verlinde algebras and indices of bundles 
over $M$ is that $\mathrm{Ind}\left(M;\cO(\tau_0-\sigma )\right)$ is the 
partition function of the surface $\Sigma$, in the 2D topological
field theory defined by $V_G(\tau _0-\sigma )$. Recall (\ref{24}) that
$V_G\otimes\bC$ is isomorphic, as an algebra, to a direct sum of copies
of $\bC$, supported on the regular Weyl orbits in $F$. The traces of 
the associated projectors are the \textit{structure constants} $\theta_f$
($f\in F^{reg}/W$) of the Frobenius algebra; their values here are 
$\left|\Delta(f)\right|^2/|F|$. The partition function of a genus 
$g$ surface is the sum $\sum \theta _f^{1-g}$, which leads to one 
version of the Verlinde formula
\begin{equation}\label{34}			
\mathrm{Ind}\left(M;\cO(\tau_0-\sigma)\right) = \sum_{f\in F^{reg}/W} 
{\left|\Delta(f)\right|^{2-2g}\cdot|F|^{g-1}}. 
\end{equation}
\begin{remark}\label{35}
The downshift by $\sigma$ stems from our use of the $\bar\partial$-index; 
the Dirac index would refer to $\cO(\tau_0)$. However, there is no definition
of the Dirac index in the singular case, and even less for the stack $\frM$.
\end{remark}

The generalization of (\ref{34}) to higher twistings is one form of 
the main conjecture. Recall the index bundle $\alpha (V)$ over $\frM$ 
associated to a representation $V$ of $G$, and call $\theta_f(t)$,
$f\in F^{reg}/W$, the structure constants of the Frobenius algebra
$\xt K_G^{\dim G}(G)\otimes \bC[[t]]$ over $\bC[[t]]$.
\begin{conjecture}\label{36}
$\mathrm{Ind}\left(\frM;\cO(\tau_0-\sigma )\otimes \exp[t\alpha (V)]\right)
= \sum\nolimits_{f\in F^{reg}/W} {\theta_f(t)^{1-g}}$.
\end{conjecture}
\begin{remark}\label{37}
\begin{trivlist}
\item(i) Expansion in $t$ allows the computation of indices of
$\cO(\tau _0-\sigma)\otimes \alpha (V)^{\otimes n}$ from (\ref{30}). 
More general expressions in the index bundles $\alpha(V)$, for 
various $V$, are easily obtained by the use of several formal parameters. 
One can also extend the discussion to include odd tautological 
generators (\ref{6}.iii), but we shall not do so here. 
\item(ii) The change in $\theta_f(t)$ is due both to the movement of the 
point $f$ under the flow, and to the change in the volume 
form, under the change of coordinates (\ref{30}).
\end{trivlist}
\end{remark}
\subsection{Transgressed twistings and the product of commutators.}\label
{38} Let us move to a more sophisticated version of the conjecture, 
which incorporates the evaluation bundles (\ref{6}.i). We refer to [FHT1], 
\S7 for more motivation, in connection to loop group representations. 
Recall the "product of commutators" map $\Pi :G^{2g}\to G$. If we remove 
a disk $\Delta$ from $\Sigma$, this map is realized by the restriction 
to the boundary of flat connections on $\Sigma \setminus\Delta$, based 
at some boundary point; the conjugation $G$-action forgets the base-point. 
The homotopy fibre of $\Pi$ is the $\Omega G$-bundle over $G^{2g}$ 
mentioned in \S\ref{3}; the actual fibre over $1\in G$ is the variety of 
based flat $G$-bundles on $\Sigma$, and its quotient by $G$-conjugation 
is $M$.

\subsection{}\label{39} The twistings $\tau$ of interest to us are 
transgressed from $G$-equivariant $B^2BU[[t]]_\otimes$ classes of a 
point. For determinantal twistings, we are looking at the transgression 
from $H^4(BG)$ to $H_G^3(G)$; and, if $G$ is simply connected, all 
determinantal twistings are so transgressed. In general, transgression 
is the integration along $S^1$ of the $B^2BU_\otimes $ class on the 
universal flat $G$-bundle over $S^1$, which is pulled back by the
classifying map of the bundle. The relevant feature of a transgressed 
twisting  $\tau$ is that its (equivariant) pull-back to $G^{2g}$, via 
$\Pi$, is trivialized by the transgression over $\Sigma -\Delta$.
This trivialization gives an isomorphism
\begin{equation}\label{40}				
K_G^*(G^{2g})\;\cong \;{}^{\Pi ^*\tau}K_G^*(G^{2g}), 
\end{equation}
\noindent which allows us to define a class ${}^\tau 1\in{}^{\Pi
^*\tau }\!K_G^0(G^{2g})$ without ambiguity. 

Recall from \S\ref{3} that the homotopy fibre of $\Pi$ over $1\in G$,
when viewed $G$-equivariantly, is represented by the stack $\frM$.
Thereon, we have two trivializations of the equivariant twisting 
$\Pi^*\tau$: one lifted from the base $\{1\}$, and one coming from 
transgression over $\Sigma \setminus\Delta$. The difference of the two 
is an element of $K^0(\frM)[[t]]_\otimes$. For the twisting $\tau_0$, 
this is the line bundle $\cO(\tau_0)$; for admissible twistings, it 
will be a formal power series in $t$, with admissible $K$-class 
coefficients. For transgressed twistings based on the exponential 
morphism (\ref{23}), we obtain the exponential $\exp[t\alpha(V)]$ of 
the index bundle $\alpha(V)$.
\begin{conjecture}\label{41}
$\mathrm{Ind}\left(\frM; \cO(\tau _0-\sigma )\otimes \exp[t\alpha (V)]\right) 
= {}^\tau Tr\left(\Pi_!{}^\tau 1\right)\in \bC[[t]]$.
\end{conjecture}
\begin{remark}\label{42}
Equality of the right-hand sides of (\ref{36}) and (\ref{41}) is part 
of the definition of the 2D field theory (Frobenius algebra) structure
on $\xt K_G^{\dim G}(G)$. The only check there is has been incorporated 
into Thm.~\ref{31}, which describes the trace map.
\end{remark}

The last formulation has the advantage of allowing us to incorporate
the evaluation bundles (\ref{6}.i). Let $W$ be another representation 
of $G$, and call $[W]$ the image class in $\xt K_G^{\dim G}(G)$.
\begin{amplification}\label{43}
$\mathrm{Ind}\left(\frM;\cO(\tau_0-\sigma)\otimes \exp[t\alpha(V)]\otimes 
E_x^*W\right) = {}^\tau Tr\left([W]\cdot \Pi_!{}^\tau 1 \right)$.
\end{amplification}

\section{The index formula by virtual localization}
In this section, we explain how the most naive localization procedure, 
from $G$ to its maximal torus, gives rise to an index formula for 
admissible $K$-classes, which agrees with Conjecture \ref{43}. 
There is an intriguing similarity here with localization methods 
used by Blau and Thomas \cite{blau} in their path-integral calculations. 
We emphasize, however, that, in twisted $K$-theory, the localization 
formula from $G$ to its maximal torus is completely rigorous, and can 
be applied to the Gysin map $\Pi_!$ of \S\ref{38} to prove the 
equivalence of Conjectures \ref{43} and \ref{51} below. The role
of the $\delta$-functions which appear in this section is played, 
on the $K$-theory side, by the skyscrapers of the sheaf $\xt K_G^
{\dim G}(G)$, when localized over the conjugacy classes. For clarity 
of the formulas we shall confine the calculation to $SU(2)$-bundles; 
no new issues appear for other simply connected groups.
\subsection{}\label{44} The maximal torus $T$ of $G=SU(2)$ is 
$S^1$, with coordinate $u$, and the moduli stack $\frM_T$ of holomorphic
$T_\bC$-bundles is $J(\Sigma)\times\bZ\times BT$, where $J(\Sigma)$ 
denotes the Jacobian variety $J_0$ of degree 0 line bundles and $BT$ 
denotes the classifying stack of $T_\bC\cong \bC^\times$. A vector bundle 
over $BT$ is a $T$-representation, and its index is the invariant 
subspace; the allowable representations, for which the index is well-defined, 
are finite-multiplicity sums of irreducibles. We have a natural 
isomorphism $H^1(J_0;\bZ) \cong H_1(\Sigma;\bZ)$, and the class $\psi 
:= \mathrm{Id}\in H^1(\Sigma)\otimes H^1(J_0)$ is the mixed $\Sigma 
\times J_0$ part of the first Chern class of the universal (Poincar\'e) 
bundle $\cP$. More precisely, denoting by $\omega$ the volume class in
$H^2(\Sigma)$, by $\eta$  the restriction of $c_1(\cP)$ to $J_0$ 
and by $\lambda \in H^2(BT)$ the Chern class of the standard 
representation of $T$, we have, for the universal bundle on
$J(\Sigma)\times\{d\}\times BT$,
\begin{equation}\label{45}	
c_1(\cP_d) = \eta + d\cdot\omega + \psi +\lambda. 
\end{equation}
\noindent Note that $e^\omega = 1 + \omega$, while $e^\psi = 1 + \psi 
- \eta\wedge\omega$, whence we get for the Chern character
\begin{equation}\label{46}			
ch(\cP_d) = (1 + \psi + d\cdot\omega -\eta\wedge\omega)\wedge e^\eta\cdot u, 
\end{equation}
\noindent having identified the Chern character of the standard representation 
with its character $u$.
\subsection{}\label{47} The ``virtual normal bundle" $\nu$ for the 
morphism $\frM_T\to \frM$ is the complex $R^*\pi_*(\mathrm{ad}\,\frg/\frt)[1]$. 
Since $\mathrm{ad}\,\frg/\frt \cong \cP^2 + \cP^{-2}$, a small 
calculation from (\ref{45}) gives, on the $d$th component
\begin{eqnarray}\label{48}			
ch\cP_d^2 &=& \left(1 + 2\psi +2 d\cdot\omega - 4\eta\wedge\omega\right)
	\cdot e^{2\eta}u^2,\nonumber \\
ch\cP_d^{-2} &=& \left(1-2\psi - 2d\cdot\omega - 4\eta\wedge\omega\right)
	\cdot e^{-2\eta }u^{-2},
\end{eqnarray}
\noindent and integrating over $\Sigma$, while remembering the shift 
by 1, gives
\begin{eqnarray}\label{49}		
ch[\nu_d] &=& u^2e^{2\eta}\cdot(g-1-2d) + u^{-2}e^{-2\eta}\cdot
	(g-1+2d)+4\eta\wedge(u^2e^{2\eta} + u^{-2}e^{-2\eta}),\nonumber \\
ch[\nu _d^*] &=& u^{-2}e^{-2\eta}\cdot(g-1-2d) + u^2e^{2\eta}\cdot
	(g-1+2d)-4\eta\wedge(u^2e^{2\eta }+u^{-2}e^{-2\eta}).
\end{eqnarray}
In the same vein, note that the Chern character of the basic line bundle
$D$, the determinant of cohomology $\det H^1\otimes\det^{-1}H^0$ of
$\cP + \cP ^{-1}$, is
\begin{equation}\label{50}				
ch(D) = e^{(2-2d)\eta}\cdot u^{-2d}. 
\end{equation}
\noindent The following conjecture describes the naive localization 
formula for the index, from $\frM$ to $\frM_T$.
\begin{conjecture}\label{51}
The index of an admissible class over $\frM$ is one-half the index of 
its restriction to $\frM_T$, divided by the equivariant $K$-theory 
Euler class of the conormal bundle $\nu^*$.
\end{conjecture}
\begin{remark}\label{52} \begin{trivlist}
\item(i) The index over $\frM_T$ is defined as integration over each
$J_d$, summation over degrees $d\in \bZ$, and, finally, selection 
of the $T$-invariant part. At the third step, we shall see that 
the character of the $T$-representation obtained from the first 
two steps is a distribution over $T$, supported at the regular 
points of $F$ (see \S\ref{24}). Miraculously, this corrects the 
problem which makes Conjecture \ref{51} impossible at first sight: 
the equivariant Euler class of $\nu ^*$ is singular at the singular 
conjugacy classes of $G$, so there is no well-defined index 
contribution over an individual component $J_d$. The sum over $d$ 
acquires a meaning by extending the resulting distribution by zero, 
on the singular conjugacy classes. 
\item(ii) The ``one-half" corrects for the double-counting of components 
in $\frM_T$, since opposite $T$-bundles induce isomorphic $G$-bundles. 
In general, we divide by the order of the Weyl group.
\end{trivlist}
\end{remark}
\subsection{}\label{53} We need the Chern character of the equivariant 
$K$-Euler class of (\ref{49}). The first two terms are sums of line 
bundles, and they contribute a multiplicative factor of
\begin{equation}\label{54}		
\left(1 - u^2e^{2\eta} \right)^{g-1+2d}\left(1 - u^{-2}e^{-2\eta }
\right)^{g-1-2d} = (-1)^{g-1}\left(ue^\eta -u^{-1}e^{-\eta }\right)
^{2g-2}\cdot \left(ue^\eta\right)^{4d}.
\end{equation} 
\noindent Now the log of the Euler class is additive, and we have
\begin{equation}\label{55}			
4\eta\cdot u^2e^{2\eta } = \left. 4\frac{d}{dx}\left({e^{x\eta}
	\cdot u^2e^{2\eta }}\right) \right|_{x=0}
\end{equation}
\noindent whence the Chern character of the Euler class of the 
remaining term is the exponential of
\begin{eqnarray}\label{56}	
-4 \frac{d}{dx}\left[\log \left(1-e^{x\eta}\cdot u^2e^{2\eta}\right)
\right.\mspace{-12mu} &+& \mspace{-12mu}\left.\log\left(1-e^{x\eta}\cdot 
u^{-2}e^{-2\eta}\right)\right]\Big|_{x=0} =
	\nonumber \\
&=& 4\frac{\eta\cdot u^2e^{2\eta}}{1 - u^2e^{2\eta}} + 4\frac{\eta\cdot
	u^{-2}e^{-2\eta }}{1 - u^{-2}e^{-2\eta}} = -4\eta,
\end{eqnarray}
\noindent and we conclude that the Chern character of the $K$-theory 
Euler class of $\nu ^*$ is
\begin{equation}\label{57}			
(-1)^{g-1}\left( {ue^\eta -u^{-1}e^{-\eta }}
\right)^{2g-2}\cdot \left( {ue^\eta } \right)^{4d}\cdot e^{-4\eta }.
\end{equation}
\subsection{}\label{58} We can now write the index formula asserted 
by Conjecture \ref{51}. For an admissible $K$-class of the form $D^
{\otimes h}\otimes\cE$, where $\cE$ is a polynomial in the classes 
(\ref{6}.i--iii), this is predicted to be the $u$-invariant part in 
the sum
\begin{equation}\label{59}			
\frac{1}{2}\,\sum_{d\in\bZ}{(-1)^{g-1}\cdot \int_{J_d} {ch(\cE)\cdot 
\frac{u^{-2(h+2)d}\cdot e^{-2(h+2)(d-1)\eta}}
{(ue^\eta - u^{-1}e^{-\eta})^{2g-2}}}}. 
\end{equation}
\noindent Noting, from (\ref{46}), that $ch(\cE)$ has a linear 
$d$-term for each factor $\alpha(V)$, it is now clear, as was 
promised in Remark (\ref{52}.i), that (\ref{59}) sums, away from
the singular points $u=\pm 1$, to a finite linear combination of 
$\delta$-functions and their derivatives, supported at the roots 
of unity of order $2(h+2)$. Integration over the torus is somewhat 
simplified by the formal change of variables $u\mapsto ue^\eta$, 
and the identification $(u-u^{-1})^2= -\left|\Delta (u)\right|^2$ 
in terms of the Weyl denominator leads to our definitive answer:
\begin{equation}\label{60}	
\mathrm{Ind}(\frM;D^{\otimes h}\otimes \cE) = \frac{1}{4\pi i}\oint {\left[ 
{\sum_{d\in\bZ} {\frac{u^{2(h+2)d}}
{\left|\Delta (u)\right|^{2g-2}}\cdot \int_{J_d}
{ch(\cE)\cdot e^{2(h+2)\eta}}}}\right]\frac{du}{u}}, 
\end{equation}
\noindent where the distribution in brackets is declared to be null 
at $u=\pm 1$.

\subsection{Example 1: Evaluation bundles.} 
Let $\cE$ be an evaluation bundle (\ref{6}.i), $\cE = E_x^*V$. 
Then,
\begin{eqnarray}\label{61}			
\int_{J_d}{ch(\cE)\cdot e^{-2(h+2)\zeta}} &=& (2h+4)^g\cdot ch_V(u),\\
\sum_{d\in\bZ} {\,u^{2(h+2)d}} &=& \frac{1}
{2h+4}\sum\limits_{\zeta ^{2h+4}=1} {\delta _\zeta (u)},
\end{eqnarray}
\noindent and the index formula (\ref{60}) reduces to the Verlinde 
formula for an evaluation bundle,
\begin{equation}\label{62}		
\mathrm{Ind}(\frM; D^{\otimes h}\otimes E_x^*V) = \sum\limits_{\scriptstyle 
{\zeta ^{2h+4}=1}\atop \scriptstyle {\Im \zeta >0}} {\frac{(2h+4)^{g-1}} 
{\left|\Delta(\zeta)\right|^{2g-2}}\cdot ch_V(\zeta)}.
\end{equation}

\subsection{Example 2: Exponentials of index bundles.} 
Let $V$ be any representation of $G$, with character a Laurent 
polynomial $f(u) = \sum f_nu^n$. We define $\dot f(u) := \sum 
n\cdot f_nu^n$, $\ddot f(u) := \sum n^2\cdot f_nu^n$. For the index
bundle $\alpha (V)$, we have
\begin{equation}\label{63}	
ch\,\alpha(V) = ch(R^*\pi_*V)-(g-1)\cdot ch(E_x^*V)
	= d\cdot \dot f(ue^\eta) - \eta\cdot \ddot f(ue^\eta ). 
\end{equation}
\noindent We compute the integral over $J_d$ for insertion in 
(\ref{60}) (after changing variables $u\mapsto ue^\zeta$):
\begin{equation}\label{64}		
\int_{\,J_d} {\exp\left[td\cdot\dot f(u)\right]\cdot \exp\left[
-\left(2h+4+t\ddot f(u)\right)\eta\right]} = \left(2h+4+t\ddot 
f(u)\right)^g\cdot \exp \left[td\cdot \dot f(u) \right], 
\end{equation}
\noindent whereupon the sum in (\ref{60}) becomes again a 
sum of $\delta$-functions
\begin{equation}\label{65}		
\sum\limits_{d\in \bZ} {\left[u\cdot\exp\left(\frac{t\dot f(u)}
	{2h+4}\right)\right]^{(2h+4)d}} = 
\sum\limits_{\zeta^{2h+4}=1}{\frac{\delta_\zeta \left(u\cdot\exp
\left(t\dot f(u)/(2h+4)\right)\right)}{2h + 4 + t \ddot f(u)}}, 
\end{equation}
\noindent the denominator arising from the change of variables inside 
the $\delta$-function. The index formula becomes now a sum over the 
solutions $\zeta_t$, with positive imaginary part, of $\zeta_t^{2h+4} 
\cdot\exp\left(t\dot f(\zeta_t)\right) = 1$:
\begin{equation}\label{66}			
\mathrm{Ind}\left(\frM; D^{\otimes h}\otimes \exp\left[t\alpha(V)
\right]\right) = \sum\limits_{\zeta_t} {\left[\frac{2h + 4 + 
t\ddot f(\zeta_t)}{\left|\Delta(\zeta_t)\right|^2}\right]^{g-1}}. 
\end{equation}
This has precisely the form predicted by Conjecture\ \ref{36}; 
note that the numerator differs from (\ref{62}) by the volume 
scaling factor in the coordinate change
\[
u\mapsto u_t := u\cdot \exp\left({t\dot f(u)}/(2h+4)\right). 
\]

\end{document}